\documentclass[12pt,notitlepage,twoside,a4paper]{amsart}
 \usepackage{amsfonts}

\usepackage{amsmath,amssymb,enumerate}

\usepackage{epsfig,fancyhdr,color}%,showkeys,amsmidx
 
\usepackage[table]{xcolor}
 
\usepackage{amssymb}
\usepackage{amsmath,amsthm}
\usepackage{latexsym}
\usepackage{amscd}
\usepackage{psfrag}
\usepackage{graphicx}
\usepackage[latin1]{inputenc}
\usepackage[all]{xy}
\usepackage[mathcal]{eucal}

\definecolor{NoteColor}{rgb}{1,0,0}

\renewcommand{\textsc}{\textcolor{red}}

\newtheorem*{theorem 1}{\rm\bf Proposition 1}
\newtheorem*{theorem 2}{\rm\bf Proposition 2}

\theoremstyle{definition}

\theoremstyle{remark}

\def\interieur#1{\mathord{\mathop{\kern 0pt #1}\limits^\circ}}

\title[Herbert Busemann (1905--1994)]{Herbert Busemann (1905--1994)
\\ A biography for his \emph{Selected Works} edition}

\author{Athanase Papadopoulos}
\address{A. Papadopoulos, Institut de Recherche Math{\'e}matique Avanc\'ee,
Universit{\'e} de Strasbourg and CNRS,
7 rue Ren\'e Descartes,
 67084 Strasbourg Cedex, France, and 
  Brown University, Mathematics Department, 
 151 Thayer Street
Providence, RI 02912, USA.}

\makeindex

 \date{}

\begin{document}

\maketitle

Herbert Busemann\footnote{Most of the information about Busemann is extracted from the following sources: 
\\
(1) An interview with Constance Reid, presumably made on April 22, 1973 and kept at the library of the G\"ottingen University. 
\\
 (2) Other documents held at the G\"ottingen University Library, published in Volume II of the present edition of Busemann's \emph{Selected Works}.  
\\
 (3) Busemann's correspondence with Richard Courant which is kept at the  Archives of New York University. A report on this correspondence, by Manfred Karbe, is published in Volume II of the present edition.
  \\
 (4) Anik\'o Sz\'abo's \emph{Vertreibung, R\"uckkehr, Wiedergutmachung} (Expulsion, return, reparation),  Wallstein Verlag, G\"ottingen, 2000, pp.~482--484. The German original as well as an English translation by M. Karbe of the article on Busemann that appeared in this book is re-edited in Volume II of the present edition. 
 
The author of the present biography is most grateful to the personnel of the Elmer Holmes Bobst Library of New York University for giving access to the correspondence between Busemann and Courant, to Manfred Karbe and Irene Zimmermann who kindly made available to him translations of this correspondence and of various other documents related to Busemann, to Hubert Goenner and the librarians of the University of G\"ottingen for providing a copy of the recording of the Reid interview and to Wallstein Verlag for the permission to re-publish the article on Busemann in the volume \emph{Vertreibung, Rückkehr, Wiedergutmachung},  as well as an English translation of it.}  was born in Berlin on May  12, 1905 and he died in Santa Ynez, County of Santa Barbara (California) on February 3, 1994, where he used to live. His first paper was published in 1930, and his last one in 1993. He wrote six books, two of which were translated into Russian in the 1960s.

Initially, Busemann was not destined for a mathematical career. His father was a very successful businessman who wanted his son to become, like him, a businessman. Thus, the young Herbert, after high school (in Frankfurt and Essen), spent  two and a half years in business. Several years later, Busemann recalls that he always wanted to study mathematics and describes this period as ``two and a half lost years of my life."

Busemann started university in 1925, at the age of 20.
 Between the years 1925 and 1930, he studied in Munich (one semester in the academic year 1925/26), Paris (the academic year 1927/28) and G\"ottingen (one semester in 1925/26, and the years 1928/1930). He also made two stays of several months in Rome, in the summer of 1930 and the winter semester of the academic year 1931/32.
 
When Busemann  enrolled in G\"ottingen University in 1926, the mathematics and physics  departments there were at their full glory. This was due  to people like David Hilbert, Emmy Noether, Bartel Leendert van der Waerden, 
Edmund Landau,
 Hermann Weyl, Gustav Herglotz, Felix Bernstein, Carl Runge, James Franck, Max Born and others. There were regular distinguished visitors at the mathematics department, like 
Pavel Alexandroff and Andre\"\i \ Kolmogorov.  Busemann's meeting with Alexandroff was crucial to him, we shall see why below. G\"ottingen had also an excellent group of German and foreign PhD students and young researchers, and by the time of Busemann,  these included William Feller, Hans Lewy, John von Neumann, Robert Oppenheimer, Werner Heisenberg, Paul Dirac, Franz Rellich, Werner Fenchel, Otto Neugebauer, Oswald Teichm\"uller, Harald Bohr, Norbert Wiener, Saunders Mac Lane, 
 and there were others.   
 
 Several years after he left G\"ottingen, Busemann shared his recollections on the teaching he received there. According to his account, Courant had so many administrative duties  that he was always unprepared for his lectures. He remembers him trying to reconstruct some proofs at the blackboard, and he found this very disturbing.  For Busemann, ``the aesthetic lecture was the main thing" and in this sense the
ideal lecturer was represented by Herglotz. Busemann adds that Alexandroff  was also a ``first rate lecturer. [...] His German was amazing. You could hear that he had an accent, but he mastered the German language perfectly." 
   Hilbert  retired in 1930,  but he continued giving a course entitled ``Introduction to philosophy on the basis of modern science." Weyl was his successor, and he lectured on differential geometry, topology and the philosophy of mathematics.

 Busemann obtained his doctoral degree on February 25, 1931.
 In a document he filled concerning his application for his PhD examination (see \cite{Goenner}), Busemann states that during his years of study in the three cities he visited, he attended courses by Pavel Alexandroff, Bessel-Hagen, Bohr, Borel, Born, 
 Carathéodory, Cartan, Courant, Franck, Graetz, Grandjot, Herglotz,
 Hilbert, Julia, Landau, Montel, Noether, Nordheim, Ostrowski, Perron,
 Walther, van der Waerden, Wegner, Weyl and Wien.
  In the last paper that Busemann wrote, \cite{BP1993}   (1993), in collaboration with Phadke, Busemann  mentions the influence of Minkowski. We read the following (p. 181):
 \begin{quote}\small
 Busemann has read the beginning of Minkowski's \emph{Geometrie der Zahlen} in 1926 which convinced him of the importance of non-Riemannian metrics. At the same time he heard a course on point set topology and learned Fr\'echet's concept of metric spaces. The older generation ridiculed the idea of using these spaces as a way to obtain results of higher differential geometry. But it turned out that a few simple axioms on distance suffice to obtain many non-trivial results of Riemannian geometry and, in addition, many which are quite inaccessible to the classical methods.
 \end{quote}
  The title of Busemann's doctoral dissertation is \emph{On
geometries in which the circles with infinite radius are the
shortest lines}, and this is also the subject of his paper \cite{Busemann1932-Uber}. The subject of the dissertation is the axiomatic characterization of Minkowski geometries (metric spaces associated to norms on finite-dimensional vector spaces). The characterization of Minkowski 
geometries is a theme that accompanied Busemann during his whole life.
In document \cite{Goenner} which we mentioned, Busemann declares that the choice of  the topic of his dissertation was inspired by the courses he followed by Alexandroff, who was then a frequent visitor in G\"ottingen. Busemann recounts the following:\footnote{All of Busemann's quotations are from the recording of the interview with Reid. The recording is corrupted at some places and we tried to reconstruct some words. We have also eliminated the questions of the interviewer and skipped some hesitations, incomplete sentences and repetitions and we have transformed some oral expressions into written English.}

\begin{quote}\small
 Officially, I took my degree  with Courant. This was only officially, in the sense that I was
 really inspired by Paul Alexandroff, who visited G\"ottingen regularly.  He gave me the idea of the subject of the thesis. I wrote it, but of course he could not be the official reviewer of my  thesis, so he was my co-referee.  
 
 I must say that my thesis was partly in protest against Courant.
  I went to Courant originally, he gave me something and it turned out to be much easier than what he thought. I did it and it became a small paper (see \cite{Busemann1930}), but
not enough for a thesis. Then he gave me something else, I did this too, 
but then Kolmogorov
came to G\"ottingen, Courant showed it to him, and Kolmogorov said: Oh that's all very fine, but it's well  known. So I became really mad, and I went away to Rome. I was angry with Courant, and wrote my thesis there, on my own.
 \end{quote}

Busemann says that he had very close contacts in Rome, and that furthermore, Courant had given him introduction letters to Levi-Civita and Enriques.  About his stay in Rome, Busemann says: ``This was the most useful thing in the sense that I decided I would never live in a fascist state. Mussolini was there since 1923, and we were in 1930."

 In a letter dated August 8, 1930, Courant informs Busemann that his ``geometric work" done in Rome is sufficient for a ``reasonable" dissertation.\footnote{An English translation of Courant's report on Busemann's dissertation is contained in the Volume II of the present \emph{Selected Works}.} The final dissertation was submitted in December 1930. After obtaining his doctorate, Busemann returned to Rome. In a letter to Lefschetz,  dated November 15, 1935, he mentions a lecture by Enriques he was attending in the winter of the academic year 1931/32, which aroused his interest in algebraic geometry.

The doctoral degree that Busemann obtained in 1931 would have led him in principle to a position of lecturer, but practically it did not. He recounts: \begin{quote}\small
Normally I would become an assistant but this was of course 
  the depression time.
Since my father had money and Courant knew it, he asked my father whether he wouldn't support me instead of having me become an assistant. This would make it possible  for some other youngster to become a mathematician. That fact did me a lot of harm  
later on,
in the sense that the German government after the Second World War tried to make good, but in my case they refused, because
I had never had a paid position. If I had a paid position I would have had a pension.
 \end{quote}
  
  Courant, as the director of the Mathematics Institute in G\"ottingen,  told Busemann that they had run over their funds  and asked him to make it possible to meet his father in order to get some financial support from him for the institute. The two men met. Busemann says: ``How they settled this situation I don't know and I thought I shouldn't ask." 
    Reid reports on this in her book on Hilbert (\cite{Reid} p. 132):
    
    \begin{quote}\small
    During the increasingly hard times, a number of informal assistantships came into being at the institute in addition to the official ones funded by the government. Often duties were vague or non-existent. Courant once gave a student a stipend because he thought the young man was on the verge of a nervous breakdown and needed a skiing vacation. He also contrived to have some students work without pay. One of them was Busemann.

  \end{quote}
 Busemann adds, concerning his former advisor:

\begin{quote}\small
Courant was rather reactionary in his mathematical outlook, and so he prevented many things which should not have been prevented. In G\"ottingen, he constantly tried to prevent the concept of Lebesgue integral. This has in the meantime conquered the whole world. He didn't see the importance of many things of modern mathematics. He had no relations at all with  algebraic geometry. [...]

The Russians had played quite a role in G\"ottingen. I believe that I was the only one who was directly inspired by them. But their course was very popular. They filled a gap. They were familiar with certain modern tendencies that were not represented in G\"ottingen, e.g. topology.
\end{quote}

On January 30, 1933, Hitler was named Chancellor of Germany. Busemann  experienced the immediate decline of the University of G\"ottingen and 
in the month of May of the same year, he left to Denmark, where he got a friendly reception by Harald Bohr, Werner Fenchel and his wife, Jakob Nielsen, the physicist and Nobel Prize winner James Frank, and others.  Busemann spent three years in Copenhagen on an unpaid position, lecturing occasionally, watching the political developments in Europe with the hope of returning to Germany after things get normal again. During these three years, Busemann does return to Germany for occasional visits, he does pay attention to what is happening there, he even intentionally  strolls  along streets in G\"ottingen, Berlin and other cities, talking to chance acquaintances to see how they think about the situation in Germany. 
When his hope of returning for good to Germany became unrealistic, Busemann decided to emigrate to the United States.  His plans to emigrate are reflected in his correspondence with Courant; cf.  \cite{Karbe1} and \cite{Karbe2}.
Courant expressed some will to help Busemann and  came up with various suggestions. For instance, he proposed that Busemann should come to the US on trial to get to know life there and see whether he feels that settling there is an option. However, there was no paid position within sight yet, and the correspondence with Courant gives the impression that the latter, in view of Busemann's financial health related to his family background, did not put particular emphasis to this point. In a letter to Courant, dated May 12, 1935, Busemann writes: ``I once thought of Baltimore, because the book of Zariski made so much impression on me." Princeton showed interest, probably because Lefschetz had the chance to encounter Busemann at some earlier event and received a good impression, but nothing resulted in a way satisfactory to Busemann. In a letter to Courant, Busemann reports that he received a letter from Oswald  Veblen in which he was invited to come to Princeton, but at the same time Veblen made clear that he  had no chance to obtain any scholarship or position whatsoever (not even an ``$epsilon$-scholarship" in the words of Busemann, letter to Courant dated November 3, 1935). 
 The colleagues in Denmark were in favor of Baltimore; only Bohr seems to prefer Princeton: ``Nobody knows Zariski." But Bohr also considers advantageously Frank's presence there. 
 In a letter to Busemann, dated October 29, 1935, Courant writes:
 \begin{quote}\small
I have now talked with both the people of Baltimore (but not directly with Zariski, who lives a very secluded life) and Lefschetz. My impression is the following: Baltimore does not seem so favorable because of  Zariski's isolation and the difficulty of having personal contact with him. Otherwise, maybe in Baltimore, the possibility of getting a scholarship for the year 1936--37  would indeed  exist. However, under the present circumstances, I would greatly prefer Princeton. Lefschetz, with whom I spoke for a long time about the case, was \emph{very} friendly, expressed himself with great warmth and appreciation, and would very much like to have you there. He suggested that you should contact him directly by letter.
\end{quote}

 During his stay in Copenhagen,  Busemann published several papers with William Feller, who had also fled Germany. Harald and Niels Bohr were actively helping German refugees in Denmark and Busemann recalls that, in Copenhagen,  ``Niels Bohr used to call up the American consulate to see how things were going for me."

Busemann's family did not leave Germany. His father stayed as one of the main directors at  Krupp, until his retirement in 1943.

In 1936, upon the recommendation of Veblen, Busemann was invited on a temporary position at Princeton's Institute of Advanced Study. Veblen was one of the first Faculty members there. Busemann spent  three years at the Institute, first as an assistant to  Marston Morse and then as a member. After that, Busemann stayed temporarily at the Johns Hopkins University and Smith College. His first  permanent position was at the Illinois Institute of Technology in Chicago, a position which he described as a ``horrible permanent job." He recalls that this was a period where  ``everybody was looking for jobs, and one had to take whatever." He says he spent ``five miserable years" 
in Chicago, from 1940 to 1945. He adds:
\begin{quote}\small
The head of the department made it difficult.
He did not like foreigners in the first place. He belonged to those people who had done a couple of good things when they were quite young
 and he  was against anyone who was too active mathematically.
On the other hand the administration forced him to take good people, and he resented them.
\end{quote}
 
 It is interesting to note that Karl Menger, who did foundational work on topology and metric geometry, joined the Faculty of 
 the Illinois Institute of Technology, in 1946 (Busemann had already left).
  
 In 1945, Busemann was appointed Assistant Professor at Smith College in Northampton.
 
 Busemann stayed in contact with Courant, and the two men had a regular correspondence, but essentially on practical matters. Talking about the institute that Courant founded later in New York, Busemann says, in his interview with Reid: ``In America, Courant tried to do again what he did before in G\"ottingen. [...] His institute is excellent but very one-sided too. The mathematics represented there goes all, or most of it, in one direction."

 In 1947,  Busemann was appointed  professor at the University of Southern California, and he spent there the rest of his career. In 1964, he was made distinguished professor. 
 
 At USC, Busemann worked in relative isolation, and practically his only collaborators were his PhD students.  His work started to be recognized only in the 1980s,  when metric geometry was revived in the West, especially by M. Gromov, and when the methods of synthetic global geometry were introduced in the study of geodesic metric spaces.   W. P. Thurston, in his approach to geometry, also started from basic principles. Before that, Busemann's work was only appreciated in the Soviet Union, where A. D. Alexandrov founded an important school on the subject, with a large number of collaborators and students. Alexandrov, like Busemann was only interested in the most basic notions of geometry. Classical problems of convexity, isoperimetry and isoepiphany became the forefront of research, and classical projective geometry took its revenge upon a certain Riemannian geometry based on linear algebra and tensor calculus in tangent spaces. In some sense, it was a return to Euclid and Archimedes.  In a  tribute to Alexandrov's memory, S. S. Kutateladze writes \cite{Kutateladze}: 
\begin{quote}\small 
Alexandroff contributed to mathematics under the slogan: ``Retreat to Euclid." He remarked that
``the pathos of contemporary mathematics is the return to Ancient Greece." 
\end{quote}
V. A. Zalgaller, in another paper dedicated to Alexandrov \cite{Zalgaller}, mentions Busemann:
 \begin{quote}\small
 In 1961, a well-known American geometer, H. Busemann (1905--1994) arrived at the IV All-Union Mathematical Congress; he also was a guest of the seminar. As a student, Busemann (almost the same age as A. D.) took part in preparation of the fundamental book by Bonnesen and Fenchel on convexity, and so he was an expert in A. D.'s area of research. He reviewed many publications of Alexander Danilovich's school for Mathematical Reviews. For the sake of that, he learned Russian and gave his talk at the Congress in Russian.
\end{quote}

In 1963, Busemann was elected foreign member of the Royal Danish Academy of Arts and Sciences.
But the most important recognition that he got for his work came from the Soviet Union, namely, the Lobachevsky Medal, which was awarded to him in 1985, ``for his innovative book 
\emph{The Geometry of Geodesics}" which he had written 30 years before. The first recipients of this prestigious prize were Sophus Lie in 1897, Wilhelm Killing in 1900 and David Hilbert in 1903.  A. D. Alexandroff  obtained it in 1951. The list of recipients also includes Weyl, Pontryagin, H. Hopf, P. S. Alexandroff,  Kolmogorov,  Hirzebruch, Arnol'd, Margulis,  Gromov and Chern.

 Busemann's work is profound. He was capable of formulating problems and working on them,  without relying on the trends that were fashionable in his time.  In an article that appeared in the Los Angeles Times on June 14, 1985, \cite{LA}, on the occasion of the attribution of the Lobachevsky prize to Busemann, the author reports the following: 
 \begin{quote}\small Few mathematicians ever make it into public consciousness, but Busemann has had a hard time even within his own field, in part, at least, because he never worked on trendy problems and never followed the crowd.
 \end{quote} The journalist quotes Busemann saying: ``If I have a
merit, it is that I am not influenced by what other people do."
   He then quotes Bob Brooks, who was Busemann's colleague at USC:  
   \begin{quote}\small 
   Tastes change a lot and interests change a lot in the space of five years. But there are people who aren't so interested in keeping up with today's fads. Busemann is
very definitely in that category.
\end{quote}
     We also read in the same article:
\begin{quote}\small
Busemann characterizes his basic mathematical approach this way: ``Any apparently difficult problem
can be done with very simple methods. This is the property of many of my things. I see a simple
geometric reason which others have overlooked."
\end{quote}
Talking about the \emph{Geometry of Geodesics} for which the Lobachevsky Medal was attributed, Busemann declares, in the  Los Angeles Times article, that the approach is more important than the results: ``The emphasis is more on the radically new approach than on the individual problem."
  
Busemann retired in 1970. In 1971, he received from USC a honorary degree of Doctor of Laws. He was a also a  linguist. He spoke several Languages, including, besides German and English, French,  Spanish, Italian, Russian, and Danish. He wrote papers or translated articles and monographs from all these languages.  Busemann published in the Mathematical Reviews a large number of  reports on articles written in Russian. He could also read Arabic, Latin, Greek and Swedish. In the article in the Los Angeles Times \cite{LA} mentioned above, the journalist quotes Busemann:  ``Every two years I read the Odyssey, I like it so much. And Plato."
Busemann was also a painter. In the same article in the Los Angeles Times, we read: ``Despite a lifelong desire to paint, Busemann never took it up, fearing that it would divert him from the arduous work of mathematics. But retirement freed him, and he built a studio in his home that is now chockablock with dozens of large canvases painted in vibrantly colored geometrical designs. His mathematical vision carried over into art."

Busemann got married in 1939. Little is known about his personal life.  In a recent letter to the author of the present article, Peter Woo, who was a student of Busemann, writes:

 \begin{quote}\small
 Busemann had a funny way of writing on the blackboard.
          Some important theorem, he would write in big letters,
from left edge to right edge, some 15 feet wide.  Then he would scribble the proof underneath, and say, ``See how easy it is,"  and then erase his proof. 
        We had to beg him not to wipe away until we finished copying it.  He would give us a puzzled look, as if we were   wasting his time.
        
        He wanted me to discover theorems about new geometrical spaces with the rule : $AB + BC \leq AC$
unless $A,B,C$ lie on a ``geodesic curve." First I found this funny. Then  I began to make conjectures, and he encouraged me to prove them, first in some particular cases, then in more and more general cases. At some point he said, ``You have done enough for the PhD. dissertation."

           He often took us to his home.
          Together with three other PhD students, we used go to his house one afternoon per week.  Each of us was supposed to present some theorem or unsolved problem, on a blackboard hanging on the wall at his backyard.  He would make some remarks on what direction we should turn to, some easy special cases we had to study first, etc.   After that, he would invite us into his house, to have a piece of pastry, and tea or coffee.  He was always positive.
          He did not tell us much about his life. He knew a lot about history of mathematics, and this tied us with the European cultural heritage.    He liked plants.  He had a cactus garden where we loved a little thing about 9" tall, 4" wide, like a ridged dark green okra, with much white hair Spreading from the top all around.  He called it ``the Old Man."  We loved it.  He had no children,  he loved us almost like his  children.  In all sense of the word, we were his disciples.
          \end{quote}

The list of Busemann's students includes 
John Beem,  John Featherstone,  Donald Glassco II,  George Lewis,   Flemming Pederson, Clinton Petty, Benson Russell,   Bhalchandra  Phadke,  J\`anos Szenthe,  Steven Weinstein, Peter Woo and Eugene Zaustinsky.

          There is not much more that can be added on Busemann's life. His work will speak for him.

\end{document}